\documentclass[12pt]{amsart}%
\usepackage{latexsym}
\usepackage{amsthm}
\usepackage{amsmath}
\usepackage{amsfonts}
\usepackage{amssymb}
\usepackage[dvips]{graphicx}
\usepackage{xypic}
\input xy
\xyoption{all}

\newtheorem{theo}{Theorem}[section]

\newtheorem{propo}[theo]{Proposition}
\newtheorem{defi}[theo]{Definition}
\newtheorem{coro}[theo]{Corollary}
\newtheorem{rem}[theo]{Remark}

\newtheorem{exam}[theo]{Example}

\newtheorem{nota}[theo]{Notation}

\newcommand\colim{\operatorname{colim}}

\newcommand\Top{\operatorname{\bf Top}}
\newcommand\PTop{\operatorname{\bf PTop}}

\newcommand\id{\operatorname{id}}
\newcommand\cof{\operatorname{cof}}

\newcommand\Mod{\operatorname{\bf Mod}}
\newcommand\Set{\operatorname{\bf Set}}

\newcommand\dSpace{\operatorname{\bf d-Space}}

\newcommand\ca{\mathcal {A}}

\newcommand\cc{\mathcal {C}}

\newcommand\ch{\mathcal {H}}
\newcommand\ci{\mathcal {I}}

\newcommand\ck{\mathcal {K}}
\newcommand\cl{\mathcal {L}}

\newcommand\crr{\mathcal {R}}

\newcommand\cB{\mathcal B}

\newcommand{\di}[1]{\vec{#1}}

\date{August 16, 2007}

\begin{document}
\title{A convenient category for directed homotopy}
\author[L. Fajstrup and J. Rosick\'{y}]
{L. Fajstrup and J. Rosick\'{y}$^*$}
\thanks{ $^*$ Supported by the Ministry of Education of the Czech republic under the projects MSM 0021622409 and
1M0545. The hospitality of the Aalborg University is gratefully acknowledged.}
\address{\newline L. Fajstrup\newline
Department of Mathematics\newline
University of Aalborg\newline
Fredrik Bajers Vej 7G, DK9220 Aalborg {\O}st, Denmark\newline
Denmark\newline
fajstrup@math.aau.dk
\newline\newline
J. Rosick\'{y}\newline
Department of Mathematics and Statistics\newline
Masaryk University, Faculty of Sciences,
\newline
Jan\'{a}\v{c}kovo n\'{a}m. 2a, 60200 Brno, Czech Republic\newline
rosicky@math.muni.cz
}
\begin{abstract}
We propose a convenient category for directed homotopy consisting of preordered topological spaces generated
by cubes. Its main advantage is that, like the category of topological spaces generated by simplices suggested
by J. H. Smith, it is locally presentable.
\end{abstract}

\maketitle

\section{Introduction}

We propose a convenient category for doing directed homotopy whose main advantage is its local presentability.
It is based on the suggestion of J. H. Smith to use $\Delta$-generated topological spaces as a convenient category
for usual homotopy. His suggestion was written down by D. Dugger \cite{D} but it turns out that it is not clear
how to prove that the resulting category is locally presentable. We will present the missing proof and, in fact,
we prove a more general result saying that for each fibre-small topological category $\ck$ and each small full
subcategory $\ci$, the category $\ck_\ci$ of $\ci$-generated objects in $\ck$ is locally presentable. In the case
of J. H. Smith, we take as $\ck$ the category $\Top$ of topological spaces and continuous maps and as $\ci$
the full subcategory consisting of simplices $\Delta_n$, $n=0,1,\dots,n,\dots$. Recall that a category $\ck$ is
topological if it is equipped with a faithful functor $U:\ck\to\Set$ to the category of sets such that one can
mimick the formation of "initially generated topological spaces" (see \cite{AHS}). The category $\dSpace$ of d-spaces 
(in the sense of \cite{G}) is topological and its full subcategory generated by suitably ordered cubes is our proposed 
convenient category for directed homotopy.

The idea of suitably generated topological spaces is quite old and goes back to \cite{W} and \cite{V} where
the aim was to get a cartesian closed replacement of $\Top$. The classical choice of $\ci$ is the category of
compact Hausdorff spaces. The insight of Smith is that the smallness of $\ci$ makes $\Top_\ci$ locally
presentable. By \cite{W} 3.3, $\Top_\Delta$ is even cartesian closed.

\section{Locally presentable categories}

A category $\ck$ is \textit{locally} $\lambda$-\textit{presentable} (where $\lambda$ is a regular cardinal) if it
is cocomplete and has a set $\ca$ of $\lambda$-presentable objects such that every object of $\ck$ is
a $\lambda$-directed colimit of objects from $\ca$. A category which is locally $\lambda$-presentable for some
regular cardinal $\lambda$ is called \textit{locally pre\-sen\-tab\-le}. Recall that an object $K$ is
$\lambda$-presentable if its hom-functor $\hom(K,-):\ck\to\Set$ preserves $\lambda$-filtered colimits. We will
say that $K$ is \textit{presentable} if it is $\lambda$-presentable for some regular cardinal $\lambda$.
A useful characterization is that a category $\ck$ is locally presentable if and only if it is cocomplete and has
a small dense full subcategory consisting of presentable objects (see \cite{AR}, 1.20).

A distinguished advantage of locally presentable categories are the following two results. Recall that, given
morphisms $f:A\to B$ and $g:C\to D$ in a category $\ck$, we write
$$
f\square g\quad\quad (f\perp g)
$$
if, in each commutative square
$$
\xymatrix{
A \ar [r]^{u} \ar [d]_{f}& C \ar [d]^{g}\\
B\ar [r]_{v}& D
}
$$
there is a (unique) diagonal $d:B\to C$ with $df=u$ and $gd=v$.

For a class $\ch$ of morphisms of $\ck$ we put
\begin{align*}
\ch^{\square}&=\{g| f\square g \mbox{ for each } f\in \ch\},\\
{}^{\square}\ch&= \{f| f\square g \mbox{ for each } g\in \ch\},\\
\ch^{\perp}&=\{g| f\perp g \mbox{ for each } f\in \ch\},\\
{}^{\perp}\ch&= \{f| f\perp g \mbox{ for each } g\in \ch\}.\\
\end{align*}
The smallest class of morphisms of $\ck$ containing isomorphisms and being closed under transfinite compositions and
pushouts of morphisms from $\ch$ is denoted as $\cof(\ch)$ while the smallest class of morphisms of $\ck$ closed under
all colimits (in the category $\ck^\to$ of morphisms of $\ck$) and containing $\ch$ is denoted as $\colim(\ch)$.

Given two classes $\cl$ and $\crr$ of morphisms of $\ck$, the pair $(\cl,\crr)$ is called a \textit{weak factorization
system} if
\begin{enumerate}
\item $\crr = \cl^\square$, $\cl = {}^\square \crr$
\end{enumerate}
and
\begin{enumerate}
\item[(2)] any morphism $h$ of $\ck$ has a factorization $h= gf$ with $f\in\cl$ and $g\in\crr$.
\end{enumerate}
The pair $(\cl,\crr)$ is called a \textit{factorization system} if condition (1) is replaced by
\begin{enumerate}
\item[(1')] $\crr = \cl^\perp$, $\cl = {}^\perp \crr$.
\end{enumerate}

While the first result below can be found in \cite{B} (or \cite{AHRT}), we are not aware of any published proof
of the second one.

\begin{theo}\label{th2.1}
Let $\ck$ be a locally presentable category and $\cc$ a set of morphisms of $\ck$. Then
$(\cof(\cc),\cc^\square)$ is a weak factorization system in $\ck$.
\end{theo}

\begin{theo}\label{th2.2}
Let $\ck$ be a locally presentable category and $\cc$ a set of morphisms of $\ck$. Then
$(\colim(\cc),\cc^\perp)$ is a factorization system in $\ck$.
\end{theo}
\begin{proof}
It is easy to see (and well known) that
$$
\colim(\cc)\subseteq{}^\perp(\cc^\perp).
$$
It is also easy to see that $g:C\to D$ belongs to $\cc^\perp$ if and only if it is orthogonal in $\ck\downarrow D$ to
each morphism $f:(A,vf)\to (B,v)$ with $f\in\cc$. By \cite{AR}, 4.4, it is equivalent to $g$ being injective to a larger
set of morphisms of $\ck\downarrow D$. Since this larger set is constructed using pushouts and pushouts in
$\ck\downarrow D$ are given by pushouts in $\ck$, $g:C\to D$ belongs to $\cc^\perp$ if and only if it is injective
in $\ck\downarrow D$ to each morphism $f:(A,vf)\to (B,v)$ with $f\in\bar{\cc}$ where $\bar{\cc}$ is given as follows
Given $f\in\cc$, we form the pushout of $f$ and $f$ and consider a unique morphism $f^\ast$ making 
the following diagram commutative
$$
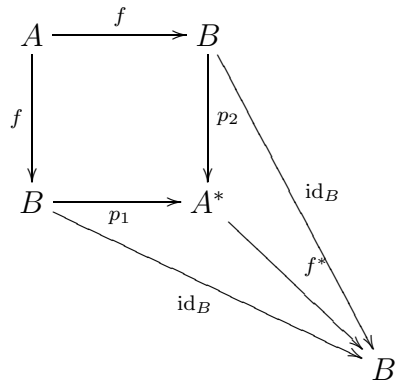
\xymatrix@C=4pc@R=4pc{
A\ar[r]^f \ar[d]_f &
B\ar[d]^{p_2}\ar[ddr]^{\id_B}&\\
B\ar[r]_{p_1}\ar[drr]_{\id_B}& A^\ast \ar[dr]^{f^\ast}&\\
&&B
}
$$
Then $f^\ast$ belongs to $\colim(\cc)$ because it is the pushout of $f:f\to\id_B$ and $f:f\to\id_B$ in $\ck^\to$ and
$f,\id_B\in\colim(\cc)$:
$$
\xymatrix@C=2pc@R=2pc{
B\ar[rrr]^{\id_B}
 \ar[ddd]_{\id_B} & & &
B\ar[ddd]^{\id_B}
\\
& A\ar[ul]_f
   \ar[r]^f
   \ar[d]_f
& B\ar[ur]_{\id_B}
   \ar[d]^{p_2} &
\\
& B \ar[dl]^{\id_B}\ar[r]_{p_1}
& A^\ast \ar[dr]^{f^\ast} &
\\
B
\ar[rrr]_{\id_B} &&& B
}
$$

Since $\bar{\cc}$ is a set, $(\cof(\bar{\cc}),\bar{\cc}^\square)$ is a weak factorization system (by \ref{th2.1}).
We have shown that
$$
\bar{\cc}^\square=\cc^\perp
$$
and
$$
\bar{\cc}\subseteq\colim(\cc).
$$
The consequence is that
$$
\cof(\bar{\cc})\subseteq\colim{\cc}.
$$
It follows from the fact that each pushout of a morphism $f$ belongs
to $\colim(\{f\})$ (see \cite{IK}, (the dual of) M13) and a transfinite
composition of morphisms belongs to their colimit closure.  In fact,
given a smooth chain of morphisms $(f_{ij}: K_i\to
K_j)_{i<j<\lambda}$ (i.e., $\lambda$ is a limit ordinal,
$f_{jk}f_{ij}=f_{ik}$ for $i<j<k$ and $f_{ij}: K_i\to K_j$ is a
colimit cocone for any limit ordinal $j<\lambda$), let $f_i:K_i \to
K$ be a colimit cocone. Then $f_0$, which is the transfinite
composition of $f_{ij}$ is a colimit in $\ck^\to$ of the chain
$$
\xymatrix@C=3pc@R=3pc{
K_0 \ar[r]^{\id_{K_0}}\ar[d]_{f_{00}}& K_0 \ar[r]
\ar[d]^{f_{01}}&\ar @{.}[r]  & K_0\ar[d]^{f_{0}}\\
K_0\ar[r]_{f_{01}}& K_1 \ar[r]& \ar@{.}[r] & K
}
$$

Thus we have
$$
\cof(\bar{\cc})\subseteq{}^\perp(\cc^\perp).
$$
Conversely
$$
{}^\perp(\cc^\perp)\subseteq{}^\square{}(\cc^\perp)={}^\square(\bar{\cc}^\square)=\cof(\bar{\cc}).
$$
We have proved that $(\colim(\cc),\cc^\perp)$ is a factorization system.
\end{proof}

\section{Generated spaces}

A functor $U:\ck\to\Set$ is called \textit{topological} if each cone
$$
(f_i:X\to UA_i)_{i\in I}
$$
in $\Set$ has a unique $U$-initial lift
$(\bar{f}_i:A\to A_i)_{i\in I}$ (see \cite{AHS}). It means that
\begin{enumerate}
\item $UA=X$ and $U\bar{f}_i=f_i$ for each $i\in I$ and
\item given $h:UB\to X$ with $f_ih=U\bar{h}_i$, $\bar{h}_i:B\to A_i$ for each $i\in I$ then $h=U\bar{h}$ for
$\bar{h}:B\to A$.
\end{enumerate}
Each topological functor is faithful and thus the pair $(\ck,U)$ is a concrete category. Such concrete categories are
called topological. The motivating example of a topological category is $\Top$.

\begin{exam}\label{ex.3.1}
{
\em
(1) A preordered set $(A,\leq)$ is a set $A$ equipped with a reflexive and transitive relation $\leq$. It means that
it satisfies the formulas
$$
(\forall x)(x\leq x)
$$
and
$$
(\forall x,y,z)(x\leq y\wedge y\leq z\to x\leq z).
$$
Morphisms of preordered sets are isotone maps, ie., maps preserving the relation $\leq$. The category of preordered
sets is topological. The $U$-initial lift of a cone $(f_i:X\to UA_i)_{i\in I}$ is given by putting $a\leq b$ on $X$
if and only if $f_i(a)\leq f_i(b)$ for each $i\in I$.

(2) An ordered set is a preordered set $(A,\leq)$ where $\leq$ is also antisymmetric, i.e., if it satisfies
$$
(\forall x,y)(x\leq y\wedge y\leq x\to x=y).
$$
The category of ordered sets is not topological because the underlying functor to sets does not preserve colimits.
}
\end{exam}

All three formulas from the example are strict universal Horn formulas and the difference between the first two and
the third one is that antisymmetry uses the equality. It was shown in \cite{R} that this situation is typical. But
one has to use the logic $L_{\infty,\infty}$ (see \cite{Di}). It means that one has a class of relation symbols whose
arities are arbitrary cardinal numbers and one uses conjunctions of an arbitrary set of formulas and quantifications
over an arbitrary set of variables. A relational universal strict Horn theory $T$ without equality then consists
of formulas
$$
(\forall x)(\varphi(x)\to\psi(x))
$$
where $x$ is a set of variables and $\varphi,\psi$ are conjunctions of atomic formulas without equality. The category
of models of a theory $T$ is denoted by $\Mod(T)$.

\begin{theo}\label{th3.2} Each fibre-small topological category $\ck$ is isomorphic (as a concrete category) to
a category of models of a relational universal strict Horn theory $T$ without equality.
\end{theo}

This result was proved in \cite{R}, 5.3. A theory $T$ can consist of a proper class of formulas. When $T$ is a set,
$\Mod(T)$ is locally presentable (see \cite{AR}, 5.30). The theory for $\Top$ is given by an ultrafilter convergence
(see \cite{R}, 5.4) and it was presented by Manes \cite{M}. This
theory is not a set of formulas. The category $\Top$ is far from being locally presentable
because it does not have a small dense full subcategory (see \cite{AR}, 1.24(7)) and no non-discrete space is
presentable (\cite{AR}, 1.14(6)).

A cone $(\bar{f}_i:A\to A_i)_{i\in I}$ is $U$-\textit{initial} if it satisfies condition (2) above. Topological
functors can be characterized as functors $U$ such that each cocone $(f_i:UA_i\to X)_{i\in I}$ has a unique $U$-final lift
$(\bar{f}_i:A_i\to A)_{i\in I}$ (see \cite{AHS}, 21.9). It means that
\begin{enumerate}
\item[(1')] $UA=X$ and $U\bar{f}_i=f_i$ for each $i\in I$ and
\item[(2')] given $h:X\to UB$ with $hf_i=U\bar{h}_i$, $\bar{h}_i:A_i\to B$ for each $i\in I$ then $h=U\bar{h}$ for
$\bar{h}:A\to B$.
\end{enumerate}
A cocone $(\bar{f}_i:A_i\to A)_{i\in I}$ is called $U$-\textit{final} if it satisfies the condition (2').

\begin{defi}\label{def3.3}
{
\em
Let $(\ck,U)$ be a topological category and $\ci$ a full subcategory of $\ck$. An object $K$ of $\ck$ is
called $\ci$-\textit{generated} if the cocone $(C\to K)_{C\in\ci}$ consisting of all morphisms from objects
of $\ci$ to $K$ is $U$-final.
}
\end{defi}
Let $\ck_\ci$ denote the full subcategory of $\ck$ consisting of $\ci$-generated objects. Using the terminology
of \cite{AHS}, $\ck_\ci$ is the \textit{final closure} of $\ci$ in $\ck$ and $\ci$ is \textit{finally dense}
in $\ck_\ci$.

\begin{rem}\label{re3.4}
{
\em
Let $\ci$ be a full subcategory of $\Top$. A topological space $X$ is $\ci$-\textit{generated} if it has
the property that a subset $S\subseteq X$ is open if and only if $f^{-1}(S)$ is open for every continuous map
$f:Z\to X$ with $Z\in\ci$. Thus we get $\ci$-generated spaces of \cite{D} in this case.

We follow the terminology of \cite{D} although it is somewhat misleading because, in the classical case of $\ci$
consisting of compact Hausdorff spaces, the resulting $\ci$-generated spaces are called $k$-spaces. A compactly
generated space should also be weakly Hausdorff (see, e.g., \cite{H}).
}
\end{rem}

\begin{propo}\label{prop3.5}
Let $(\ck,U)$ be a topological category and $\ci$ a full subcategory. Then $\ck_\ci$ is coreflective in $\ck$ and
contains $\ci$ as a dense subcategory.
\end{propo}
\begin{proof}
By \cite{AHS}, 21.31, $\ck_\ci$ is coreflective in $\ck$. Since $\ci$ is finally dense in $\ck_\ci$, it is
dense.
\end{proof}

The coreflector $R:\ck\to\ck_\ci$ assigns to $K$ the smallest $\ci$-generated object on $UK$.

A concrete category $(\ck,U)$ is called \textit{fibre-small} provided that, for each set $X$, there is only a set
of objects $K$ in $\ck$ with $UK=X$.

\begin{theo}\label{th3.6}
Let $(\ck,U)$ be a fibre-small topological category and let $\ci$ be a full small subcategory of $\ck$. Then the category
$\ck_\ci$ is locally presentable.
\end{theo}
\begin{proof}
By \ref{th3.2}, $\ck$ is concretely isomorphic to $\Mod(T)$ where $T$ is a relational universal strict Horn theory without
equality. We can express $T$ as a union of an increasing chain
$$
T_0\subseteq T_1\subseteq\dots T_i\subseteq\dots
$$
of subsets $T_i$ indexed by all ordinals. The inclusions $T_i\subseteq T_j$, $i\leq j$ induce functors
$H_{ij}:\Mod(T_j)\to\Mod(T_i)$ given by reducts. Analogously, we get functors $H_i:\Mod(T)\to\Mod(T_i)$ for each $i$.
All these functors are concrete (i.e., preserve underlying sets) and have left adjoints
$$
F_{ij}:\Mod(T_i)\to\Mod(T_j)
$$
and
$$
F_i:\Mod(T_i)\to\Mod(T).
$$
These left adjoints are also concrete and $F_i(A)$ is given by the $U$-initial lift of the cone
$$
f:U_i(A)\to U(B)
$$
consisting of all maps $f$ such that $f:A\to H_i(B)$ is a morphism in $\Mod(T_i)$. The functors $F_{ij}$ are
given in the same way. Since these left adjoints are concrete, they are faithfull and it immediately follows
from their construction that they are also full. Thus we have expressed $\Mod(T)$ as a union of an increasing
chain of full coreflective subcategories
$$
\Mod(T_0)\subseteq\Mod(T_1)\subseteq\dots\Mod(T_i)\subseteq\dots
$$
indexed by all ordinals. Moreover, all these coreflective subcategories are locally presentable.

Let $\ci$ be a full small subcategory of $\ck$. Then there is an ordinal $i$ such that $\ci\subseteq\Mod(T_i)$.
Consequently, $\ck_\ci\subseteq\Mod(T_i)$ and thus $\ck_\ci$ is a full coreflective subcategory of
a locally presentable $\Mod(T_i)$ having a small dense full subcategory $\ci$. Since $\ci$ is a set, there is
a regular cardinal $\lambda$ such that all objects from $\ci$ are $\lambda$-presentable in $\Mod(T_i)$ (see
\cite{AR}, 1.16). Since $\ck_\ci$ is closed under colimits in $\Mod(T_i)$, each object from $\ci$ is
$\lambda$-presentable in $\ck_\ci$. Hence $\ck_\ci$ is locally $\lambda$-presentable.
\end{proof}

\begin{coro}\label{cor3.7}
Let $\ci$ be a small full subcategory of $\Top$. Then the category $\Top_\ci$ is locally presentable.
\end{coro}

\begin{rem}\label{re3.8}
{
\em
Let $\ck$ be a category such that the coreflective closure $\ck_\ci$ of each small full subcategory $\ci$ of $\ck$
is locally presentable. Then $\ck$ is a union of a chain
$$
\ck_0\subseteq\ck_1\subseteq\dots\ck_i\subseteq
$$
of full coreflective subcategories which are locally presentable. It suffices to express $\ck$ as a union of a chain
$$
\ci_0\subseteq\ci_1\subseteq\dots\ci_i\subseteq
$$
of small full subcategories and pass to
$$
\ck_{\ci_0}\subseteq\ck_{\ci_1}\subseteq\dots\ck_{\ci_i}\subseteq
$$
}
\end{rem}

\begin{theo}\label{th3.9}
Let $\ci$ be a full subcategory of $\Top$ containing discs $D_n$ and spheres $S_n$, $n=0,1,\dots$. Then
the category $\Top_\ci$ admits a cofibrantly generated model structure, where cofibrations and weak equivalences
are the same as in $\Top$.
\end{theo}
\begin{proof}
Analogous to \cite{H}, 2.4.23.
\end{proof}

\section{Generated ordered spaces}

In order to get our convenient category for directed homotopy, we have to replace $\Top$ by a suitable category
of ordered topological spaces. We have considered two such
categories:
\begin{itemize}
\item The category $\PTop$ of preordered topological spaces. Its
objects are topological spaces whose underlying set is
preordered. Morphisms are continuous maps $f$ s.t. $x\leq y\Rightarrow
f(x)\leq f(y)$. 

\item The category $\dSpace$ of topological spaces $X$
  with a set of paths $\vec{P}(X)\subset X^I$ (see \ref{def4.1}).
\end{itemize}
 These are all topological categories, i.e., the forgetful functor to
 $\Set$ is topological, and they are directed. We would like to
 have directed loops in the category, i.e., the circle $S^1$ with
 counterclockwise direction. In $\PTop$ we require transitivity, and
 hence, a relation relating pairs of points on the circle $e^{i\theta}\leq
 e^{i\phi}$ when $\theta\leq \phi$, will be the trivial relation in $\PTop$ 

In $\dSpace$, \cite{G} the directions are represented in the allowed paths and
not as a relation on the
space itself. On a d-space, $(X,\vec{P}(X))$ the relation $x\leq y$ if
there is $\gamma\in \vec{P}(X)$ s.t. $\gamma(0)=x$ and $\gamma(1)=y$
is  gives a functor from $\dSpace$ to
$\PTop$. In the other direction, the increasing continuous maps from
$\vec{I}$ to a space in $\PTop$ will give a set of
dipaths, hence a functor to $\dSpace$. 

\begin{defi}\label{def4.1}
{
\em
The objects in $\dSpace$ are pairs $(X,\vec{P}(X))$, where $X$ is a
topological space and
$\vec{P}(X)\subset X^I$ satisfies
\begin{itemize}
\item All constant paths are in $\vec{P}(X)$
\item $\vec{P}(X)$ is closed under concatenation and increasing
  reparametrization.
\end{itemize}
$\vec{P}(X)$ is called the set of dipaths or directed paths.

A morphism $f:(X,\vec{P}(X))\to (Y,\vec{P}(Y))$ is a continuous map
$f:X\to Y$ s.t. $\gamma\in\vec{P}(X)$ implies $f\circ\gamma\in\vec{P}(Y)$
}
\end{defi}

In $\dSpace$ we do have directed circles.

\begin{theo}\label{th4.2}
$\dSpace$ is a topological category.
\end{theo}
\begin{proof}
Let $T$ be a relational universal strict Horn theory without equality giving $\Top$
and using relation symbols $R_j$, $j\in J$. We add a new continuum-ary relation
symbol $R$ whose interpretation is the set of directed paths. We add to $T$ the following
axioms:
\begin{enumerate}
\item[(1)] $(\forall x)R(x)$ where $x$ is the constant,
\item[(2)] $(\forall x,y,z)(\bigwedge\limits_{0<i\leq\frac{1}{2}} z_t=x_t
\wedge\bigwedge\limits_{0<i\leq\frac{1}{2}} z_{\frac{1}{2}+i}=y_i
\wedge x_1= y_0\wedge R(x)\wedge R(y)\to R(z)),$
\item[(3)] $(\forall x)(R(x)\to R(xt)$ where $t$ is an increasing reparametrization,
\item[(4)] $(\forall x) (R(x)\to R_j(xa))$ where $j\in J$ and $I$ satisfies $R_j$ for $a$.
\end{enumerate}
The resulting relational universal strict Horn theory axiomatizes d-spaces. In fact, (1) makes
each constant path directed, (2) says that directed paths are closed under concatenation, (3)
says that they are closed under increasing reparametrization and (4) says that they are continuous.
\end{proof}

\begin{rem}\label{re4.3}
{
\em
(i) A d-space is called \textit{saturated} if it satisfies the converse implication to (3):

(5) $(\forall x)(R(xt)\to R(x)$ where $t$ is an increasing reparametrization

\noindent
It means that a path is directed whenever some of its increasing re\-pa\-ra\-met\-ri\-za\-ti\-ons 
is directed. Thus saturated d-spaces also form a topological category.

(ii) There is, of course, a direct proof of \ref{th4.2}. By \cite{AHS}, 21.9, it suffices to see that 
the forgetful functor $U:\dSpace \to \Set$ satisfies:
For any cocone $(f_i:UA_i\to X)$ there is a unique \emph{$U$-final lift}
$(\bar{f}_i:A_i\to A)$, i.e., there is a unique $\dSpace$ structure on
$X$ such that $h:X\to UB$ is a d-morphism whenever $h\circ f_i$ is a
d-morphism for all $i$. The topology is defined by $V$ open
if and only if $f_i^{-1}(V)$ open for all $i$. Let $\vec{P}(A)$ be the
closure under concatenation and increasing reparametrization of the set of
all constant paths and all $f_i\circ\gamma$ where
$\gamma\in\vec{P}(A_i)$. It is not hard to see, that this is a
$U$-final lift.
}
\end{rem}

\begin{coro}\label{cor4.4}
Let $\ci$ be a small full subcategory of $\dSpace$. Then the category $\dSpace_\ci$ is locally presentable.
\end{coro}

\begin{defi}\label{def4.5}
{
\em
Let $\cB$ be the full subcategory of $\dSpace$ with objects all cubes
$I_1\times I_2\times \ldots \times I_n$ where $I_k$ is either the unit
interval with the trivial order (i.e., $a\leq b$ for all $a,b$) or the unit interval with the standard
order. The (pre)order on $I_1\times I_2\times \ldots \times I_n$ is the
product relation. The dipaths are the increasing paths wrt. this relation.
}
\end{defi}

\begin{nota}\label{not4.6}
{
\em
Let $I$ denote the unit interval with the trivial order and let $\di{I}$ denote the unit interval with the standard order.
}
\end{nota}

\begin{coro}\label{cor4.7}
The category ${\dSpace}_{\cB}$ is locally presentable.
\end{coro}

We consider the category ${\dSpace}_{\cB}$ a suitable framework for
studying the directed topology problems arising in concurrency. One
reason for this is, that the geometric realization of a cubical
complex is in ${\dSpace}_{\cB}$. These are geometric models of Higher
Dimensional Automata, see \cite{FGR}. In \cite{FGR}, the directions on
the spaces are given via a {\em{local partial order}} and not as d-spaces, 
but the increasing paths wrt. the local partial order are
precisely the dipaths in the d-space structure. 

For directed homotopy theory, this category is also suitable:

\begin{defi}\label{def4.8}
{
\em
Let $f,g:X\to Y$ be d-maps. A d-homotopy \cite{G} is a d-map
$H:X\times \vec{I}\to Y$ s.t. $H(x,0)=f(x)$ and $H(x,1)=g(x)$; the
d-homotopy equivalence relation is the reflexive transitive hull of
this relation. A
d-homotopy of dipaths $\gamma$, $\mu$ with common initial and final
points is a d-map $H:\vec{I}\times \vec{I}\to Y$
s.t. $H(t,0)=\gamma(t)$, $H(t,1)=\mu(t)$ and $H(0,s)=\gamma(0)=\mu(0)$
and $H(1,s)=\gamma(1)=\mu(1)$. 

A dihomotopy \cite{FGR} is unordered along the homotopy coordinate: $H:X\times
I\to Y$. This gives an equivalence relation without closing
off. Dihomotopic dipaths are defined as above - with fixed endpoints.
}
\end{defi}

Since we allowed both the trivially ordered interval and the naturally
ordered interval in $\cB$, the category $\dSpace_{\cB}$ is convenient
for both kinds of directed homotopy.
 
Globes have been considered as models for higher dimensional automata,
in \cite{GG}. A globe on a non-empty (d-)space $X$ is the unreduced suspension
$X\times \vec{I}/(x,1)\sim *_1, (x,0)\sim *_0$. If $X$ is in $\dSpace_{\cB}$ then clearly so is the globe of
$X$ as a coequalizer. The globe of the empty set is the d-space of two
disjoint points, which is also in $\dSpace_{\cB}$

The elementary globes, the globe of an unordered ball, are
equivalent to the globe of an unordered cube, which is in our category.

\section{Dicoverings}
In \cite{F1}, dicoverings, i.e., coverings of directed topological spaces are introduced
as a counterpart of coverings in the undirected case. The categorical
framework there is (subcategories of) locally partially ordered spaces. It turns out,
that it is not obvious which category, one should choose to get
universal dicoverings. With the
framework here, we have a setting which on the one hand is much more
general than the almost combinatorial one of cubical sets, and on the
other hand, it is not as general as locally partially ordered
topological spaces, where dicovering theory is certainly not well
behaved.
 In \cite{F1} we consider dicoverings with respect to a basepoint, a
 fixed initial point.

\begin{defi}\label{def5.1}
{
\em
Let $p:Y\to X$ be a morphism in $\dSpace$, let $x_0\in
  X$. Then $p$ is a dicovering wrt. $x_0$ if for all $y_0\in
  p^{-1}(x_0)$ and all $\gamma\in\vec{P}(X)$ with $\gamma(0)=x_0$,
  there is a unique lift $\hat{\gamma}$ with $\hat{\gamma}(0)=y_0$:
$$\xymatrix{{\{0\}}\ar[r]\ar@{^{(}->}[d]&Y\ar[d]^{p}\\\di{I}\ar@{-->}[ur]^{\hat{\gamma}}\ar[r]^\gamma
  &X}$$
And for all $H:I\times\vec{I}\to X$ with $H(s,0)=x_0$, there is a  unique lift $\hat{H}$:
$$\xymatrix{(I\times\{0\},I\times\{0\})\ar[r]\ar@{^{(}->}[d]&(Y,y_0)\ar[d]^{p}\\(I\times\di{I},I\times\{0\})\ar@{-->}[ur]^{\hat{H}}\ar[r]^H
  &(X,x_0)}$$
}
\end{defi}

In the present framework however, we will consider lifting properties
wrt. all initial points:
Let $J$ be the coequalizer
$$\xymatrix{I\ar@/^/[r]^g\ar@/_/[r]_f &I\times\vec{I}\ar[r] & J}$$
where $f(x)=(0,0)$ and $g(x)=(x,0).$

\begin{defi}\label{def5.2}
{
\em
Let $p:Y\to X$ be a morphism in $\dSpace$. Then $p$ is a
  dicovering, if for all $\gamma\in\vec{P}(X)$ there is a unique lift
  $\hat{\gamma}$
$$\xymatrix{{\{0\}}\ar[r]\ar@{^{(}->}[d]&Y\ar[d]^{p}\\\di{I}\ar@{-->}[ur]^{\hat{\gamma}}\ar[r]^\gamma
  &X}$$

and for all $H:J\to X$ there is a unique lift
$$\xymatrix{{*}\ar[r]\ar@{^{(}->}[d]&Y\ar[d]^{p}\\ J\ar@{-->}[ur]^{\hat{\gamma}}\ar[r]^\gamma
  &X}$$
where $*$ is the point $ (x,0)\in J.$
}
\end{defi}

Hence, a dicovering is a morphism $p:Y\to X$ which has the unique right lifting property
with respect to the inclusions $\mathcal{C}=\{0\to \vec{I},
*\to J\}$. Hence

\begin{propo}\label{prop5.3} 
A morphism $p:Y\to X$ in $\dSpace$ is a dicovering if and only if it is in $\mathcal{C}^{\perp}.$
\end{propo}

\begin{defi}\label{def5.4}
{
\em
A \textit{universal dicovering} of $X\in \dSpace_{\cB}$ is a morphism $\pi:\tilde{X}\to X$ such that for any dicovering
$p:Y\to X$ in $\dSpace_{\cB}$, there is a unique morphism $\phi:\tilde{X}\to Y$ such that $\pi=p\circ\phi.$
}
\end{defi}

\begin{coro}\label{cor5.5}
Let $X\in \PTop_{\cB}$. Then there is a universal
  dicovering $\pi:\tilde{X}\to X$, and it is unique.
\end{coro}
\begin{proof} This follows from \ref{th2.2}, since $\dSpace_{\cB}$ is
locally presentable. 
Let
$$
0 \xrightarrow{\ w \ } \tilde{X}\xrightarrow{\ u \ } X
$$
be the $(\colim(\mathcal C),\mathcal C^\bot)$ factorization of the unique morphism from the initial object
$0$ (the empty set) to $X$. Then $u:\tilde{X}\to X$ is a universal dicovering of $X$. In fact, each dicovering 
$$
v:Y\to X
$$
has a unique factorization through $u$. It suffices to apply the unique right lifting property to
$$
\xymatrix@C=3pc@R=3pc{
\tilde{X} \ar[r]^{u} & X \\
0 \ar [u]^{} \ar [r]_{} &
Y \ar[u]_{v}
}
$$
\end{proof}

In \cite{F1}, we construct a ``universal'' dicovering $\pi:\tilde{X}\to X$
by endowing the set
of dihomotopy classes of dipaths initiating in a fixed point $x_0$
with a topology and a local partial order. If all points in $X$ are
reachable by a directed path from $x_0$ and if $X$ is in $\dSpace_{\cB}$
the construction here and the underlying d-space of the locally
partially ordered space $\tilde{X}$ in
\cite{F1} should coincide, but we do not have a proof of this yet.

\end{document}